\documentclass[sigconf,natbib=false]{acmart}

\usepackage[american]{babel}
\usepackage[T1]{fontenc}
\usepackage[utf8]{inputenc}
\usepackage[american]{isodate}
\usepackage{soul}

\usepackage[ruled,norelsize]{algorithm2e}

\usepackage{amsmath}
\usepackage{amsthm}
\usepackage{bm}

\usepackage{graphicx}
\usepackage{pgf, tikz}
\usepackage{pgfplots}
\pgfplotsset{compat=1.17}
\usepgfplotslibrary{statistics}
\usepackage{placeins}

\usepackage{subfig}
\usepackage{pgfplots}

\usepackage[datamodel=acmdatamodel,style=acmnumeric]{biblatex}
\usepackage[autostyle=true]{csquotes}
\usepackage{listingsutf8}
\usepackage{multirow}
\usepackage{orcidlink}
\usepackage{seqsplit}
\usepackage{siunitx}
\usepackage{xcolor}

\definecolor{sign}{HTML}{b02a2d}
\definecolor{direction}{HTML}{007900}
\definecolor{regime}{HTML}{8c399e}
\definecolor{characteristic}{HTML}{1f5dc2}
\definecolor{mantissa}{HTML}{636363}
\definecolor{error}{HTML}{BD002A}
\definecolor{cellbg}{HTML}{EDEDED}

\definecolor{p-sign}{HTML}{FF5454}
\definecolor{p-regime}{HTML}{CC9966}
\definecolor{p-regime-term}{HTML}{996633}
\definecolor{p-exponent}{HTML}{0080FF}
\definecolor{p-fraction}{HTML}{000000}

\definecolor{posit}{HTML}{b02a2d}
\definecolor{bfloat16}{HTML}{007900}
\definecolor{takum}{HTML}{1f5dc2}
\definecolor{float}{HTML}{636363}

\makeatletter
\def\lst@makecaption{%
  \def\@captype{table}%
  \@makecaption
}
\makeatother

\setcopyright{acmlicensed}
\copyrightyear{2025}
\acmYear{2025}
\acmDOI{XXXXXXX.XXXXXXX}

\acmConference[SC '25]{International Conference for High Performance Computing, Networking, Storage and Analysis}{November 16--21,
  2025}{St. Louis, MO}

\acmBooktitle{SC '25: Proceedings of the International Conference for High Performance Computing, Networking, Storage and Analysis, November 16--21, 2025, St. Louis, MO}
\acmISBN{978-1-4503-XXXX-X/2025/11}

\begin{document}

\title{Spectral Methods via FFTs in Emerging Machine Number Formats: OFP8, Bfloat16, Posit, and Takum Arithmetics}

\author{Laslo Hunhold}
\email{hunhold@uni-koeln.de}
\orcid{0000-0001-8059-0298}
\affiliation{%
  \institution{University of Cologne}
  \city{Cologne}
  \country{Germany}
}

\author{John Gustafson}
\email{jlgusta6@asu.edu}
\orcid{0000-0002-2957-1304}
\affiliation{%
  \institution{Arizona State University}
  \city{Tempe}
  \state{Arizona}
  \country{USA}
}

\begin{abstract}
The Fast Fourier Transform (FFT) is one of the most widely used algorithms in high performance computing, with critical applications in spectral analysis for both signal processing and the numerical solution of partial differential equations (PDEs). These data-intensive workloads are primarily constrained by the memory wall, motivating the exploration of emerging number formats---such as OFP8 (E4M3 and E5M2), bfloat16, and the tapered-precision posit and takum formats---as potential alternatives to conventional IEEE 754 floating-point representations.
\par
This paper evaluates the accuracy and stability of FFT-based computations across a range of formats, from 8 to 64 bits. Round-trip FFT is applied to a diverse set of images, and short-time Fourier transform (STFT) to audio signals. The results confirm posit arithmetic's strong performance at low precision, with takum following closely behind. Posits show stability issues at higher precisions, while OFP8 formats are unsuitable and bfloat16 underperforms compared to float16 and takum.
\end{abstract}




\maketitle

\section{Introduction}\label{sec:introduction}
The Fast Fourier Transform (FFT) is a foundational algorithm in high-performance computing. Its importance extends beyond traditional applications in signal analysis, also, for instance, playing a central role in spectral methods for solving partial differential equations (PDEs). As processor performance continues to outpace improvements in memory bandwidth---a phenomenon known as the \enquote{memory wall}---there is increasing pressure to adopt lower-precision arithmetic in order to reduce memory traffic and improve computational efficiency. However, standard IEEE 754 floating-point formats are known to perform poorly in low-precision settings, particularly with respect to dynamic range and rounding behavior.
\par
In response, a variety of emerging number formats has been proposed, offering the potential to improve both accuracy and performance at reduced bit widths. Among these are the 8-bit E4M3 and E5M2 types defined in the Open Compute Project’s 8-bit Floating Point (OFP8) specification \cite{ofp8}. E4M3 uses four exponent bits and three fraction bits, while E5M2 allocates five bits to the exponent and two to the fraction, trading precision for greater dynamic range. Both formats have been adopted in Intel’s latest vector instruction set architecture, AVX10.2 \cite{intel-avx10.2}. Also supported in AVX10.2 is the 16-bit \texttt{bfloat16} format \cite{bfloat16}, developed by Google Brain. Widely implemented in hardware accelerators, \texttt{bfloat16} has become a standard format for representing weights in deep learning applications.
\par
Departing from IEEE 754-based designs, posit arithmetic, introduced by \textsc{Gustafson} et al. in 2017 \cite{posits-beating_floating-point-2017, posits-standard-2022}, employs a tapered precision approach. This scheme assigns greater precision to values near one, using a variable-length encoding of the exponent. In addition to improved precision characteristics, this design offers advantages in numerical behavior and hardware design \cite{posit-na-1, posit-na-2, posit-hardware_cost-2019}. Extending this concept, takum arithmetic \cite{2024-takum} introduces an alternative tapered-precision format that offers significantly increased dynamic range, particularly at very low bit widths. It achieves higher precision for both small and large magnitudes, at the expense of slightly reduced precision near one, while simplifying hardware implementation through a modified exponent encoding scheme. In this work, we focus exclusively on the floating-point variant---known as linear takums---and, for brevity, refer to them simply as takum throughout the paper.
\par
This paper evaluates the numerical performance of these number formats through a series of FFT-based experiments. Our benchmarks include standard spectral-method solutions to PDEs, round-trip FFT on large image datasets, and short-time Fourier transforms (STFT) applied to audio signals. To ensure practical relevance, we do not tailor algorithms or inputs to accommodate low-precision arithmetic. Most of the examined formats are implemented in software, and we do not evaluate computation time or power consumption; our focus is exclusively on numerical behavior.
\par
To date, no comprehensive study has compared all of these formats side-by-side on representative datasets and algorithms. While applications of \texttt{bfloat16} have been explored in specific domains such as astronomy \cite{bfloat16-astronomy}, and posit arithmetic has been studied in limited low-precision settings using \enquote{quire} accumulators for exact dot products \cite{posit-lossless_ffts}, such approaches preclude a uniform comparison. Other studies involving posit arithmetic rely on synthetic datasets and small-scale validation benchmarks \cite{posit-spectral_analysis}. This work provides the first broad evaluation of emerging arithmetic formats across diverse FFT-based applications using consistent test conditions.
\par
The remainder of this paper is organized as follows. Section~\ref{sec:methods} describes the setup and methods used to benchmark the formats. Section~\ref{sec:results} outlines the main results, including detailed analyses and visualizations. Finally, Section~\ref{sec:conclusion} summarizes our findings and offers conclusions.
\section{Experimental Methods}\label{sec:methods}
This work evaluates the numerical performance of FFTs across three domains: spectral solvers, demonstrated using two canonical test problems (the 1D heat equation and the \textsc{Poisson} equation); round-trip 2D FFTs applied to image data; and short-time Fourier transforms (STFT) applied to audio signals.
\par
To support these evaluations, a benchmark framework named NGAFFTB (Next-Generation Arithmetic Fast Fourier Transform Benchmarks) is presented in this work. Designed for scalability and extensibility, NGAFFTB provides a structured basis for systematic comparisons across number formats. The framework comprises four core components: a dataset preparer for automatic downloading and preprocessing of image and audio datasets, a unified interface for executing various experiment types, an implementation of the PDE test problems, and implementations of the image round-trip FFT and audio STFT experiments under investigation. The following sections provide a detailed description of each component.
\subsection{Image and Audio Datasets}
The choice of a dataset is a difficult task, as one can get very specific and distinguish between classes of image and audio types. However, considering a high-performance computing (HPC) context with diverse input data, it is reasonable to choose large, diverse datasets and evaluate the overall performance. While one may argue that knowing which arithmetic performs better for a specific dataset type could be useful---perhaps apart from AI applications---it is in general unwarranted to buy or build specialised hardware in an HPC context that is tailored to a single application.
\par
The recent rise of AI has fortunately led to the creation of many large, diverse datasets for training and validation. However, given that we still use mostly software-emulated arithmetic, the dataset should not be excessively large either.
Specifically, in this study we choose the image dataset DIV2K~\cite{div2k}, a diverse full-colour image dataset comprised of 900 high-resolution images (800 for training, 100 for validation). We use the officially provided factor-two bicubic downscaled version to reach a reasonable per-image resolution of 1020 by 678 pixels. In the NGAFFTB toolbox, the respective dataset archives are downloaded and extracted to the folder \texttt{data/image/}~\cite[src/generate\_image\_dataset.sh]{ngafftb}.
\par
For the audio dataset, we use ESC-50~\cite{esc-50}, a collection of 2000 environmental audio recordings intended for benchmarking sound classification methods. It consists of 5-second long (WAV, \SI{44.1}{\kilo\hertz}, mono) recordings spanning 50 semantic classes (40 examples per class), grouped into five broad categories: animals, natural soundscapes and water sounds, human non-speech sounds, interior/domestic sounds, and exterior/urban noises. Given their nature, these sounds exhibit highly diverse frequency spectra and are an ideal choice for benchmarking FFT algorithms, while also being representative of real-world, diverse audio input. In the NGAFFTB toolbox, the respective dataset archives are downloaded and extracted to the folder \texttt{data/audio/}~\cite[src/generate\_audio\_dataset.sh]{ngafftb}.
\subsection{Common Experiment Interface}
The overall approach is consistent across all experiments, each of which consists of multiple runs. In each experiment run, a set of parameters is used to prepare common data (e.g., reference solutions) for all arithmetic data types. The experiment is then executed in parallel across all target arithmetic types, and the output is compared to the reference solution to yield error metrics \cite[src/Experiments.jl]{ngafftb}.
\par
For each numeric format under evaluation, the input data, including parameters, is fully converted to the corresponding format. Naturally, the specifics of this conversion vary depending on the experiment and will be detailed in later sections. All computations are carried out entirely in the target data type.
\par
The output, which is in the target data type, is subsequently converted to \texttt{float128} to allow for meaningful comparison with the precomputed reference solution. The construction of these reference solutions will be discussed in the respective sections.
\subsection{Spectral PDE Experiments}
This work evaluates two partial differential equation (PDE) test problems implemented in the toolbox \cite[src/PDE.jl]{ngafftb} used in the experiments: the one-dimensional heat equation (parabolic PDE) and the two-dimensional \textsc{Poisson} equation (elliptic PDE). The heat equation is defined for \( u \colon [0,1]^2 \to \mathbb{R} \) and thermal diffusivity \( \alpha > 0 \) as
\begin{equation}
	\frac{\partial u(t,x)}{\partial t} = \alpha \cdot  \Delta u(t,x) :=
	\alpha \cdot \frac{\partial^2 u(t,x)}{\partial x^2}.
\end{equation}
As initial condition, a \textsc{Gauss}ian pulse centered at 0.5 is used:
\begin{equation}
	u(0,x) = \exp\!\left(\frac{{(x - 0.5)}^2}{2 \sigma^2}\right),
\end{equation}
a common choice. The parameter $\sigma $ is set to 0.1. Discretization is introduced through two parameters, $N_t, N_x \in \mathbb{N}_1$, splitting the time domain $[0,1]$ into $N_t + 1$ segments and the spatial domain into $N_x + 1$ segments.
\par
The \textsc{Poisson} problem is defined over $[0,1]^2$ with Dirichlet boundary conditions. Given functions $u, f \colon [0,1]^2 \to \mathbb{R}$, the equation is specified as
\begin{equation}
	\begin{cases}
		-\Delta u(x,y) = f(x,y) & (x,y) \in {(0,1)}^2\\
		u(x,y)=0 & (x,y) \in \partial {[0,1]}^2.
	\end{cases}
\end{equation}
The right-hand side is constructed to yield a \textsc{Gauss}ian pulse in the solution. Let
\begin{equation}
	r(x,y) := \sqrt{{(x - 0.5)}^2 + {(x - 0.5)}^2},
\end{equation}
then the source term is defined as
\begin{equation}
	f(x,y) := \exp\!\left(
		\frac{-r^2(x,y)}{2 \sigma^2}
	\right)
	\frac{r^2(x,y) - 2\sigma^2}{\sigma^4},
\end{equation}
with $\sigma > 0$.
The domain is discretized into a regular grid by dividing each spatial dimension into $N_x + 1$ segments.
\par
While it may be tempting to compare results against analytical solutions to quantify numerical error, such comparisons are confounded by the inherent discretization error introduced by finite $N_x$ and $N_t$, which persists even under exact arithmetic. Therefore, this work employs a high-precision reference solution computed using \texttt{float128} to evaluate the numerical accuracy of each arithmetic format under consideration~\cite[src/Experiments.jl]{ngafftb}. 
The explicit experiment setups, including parameter specifications and invocations of the common experiment interface, are located in \texttt{src/solve\_pde-*.jl}.
\subsection{FFT Experiments}
The FFT experiments are structured as follows: For the image dataset, each input image is loaded into a matrix of RGB samples. These samples are separated into their respective channels and converted to the numeric format under evaluation. An FFT is applied independently to each channel, followed by an inverse FFT. The reconstructed channels are then recombined into a full-color image~\cite[src/Experiments.jl, line 387]{ngafftb}. The reference \enquote{solution} is the original input image, with each channel cast to \texttt{float128}.
\par
For the audio dataset, a different approach is taken, reflecting the non-stationary nature of audio signals. A common technique for analyzing such signals is the Short-Time Fourier Transform (STFT), in which the audio sample is divided into short, approximately stationary segments of fixed length, known as the window size. These segments overlap according to a specified hop size; for instance, a hop size equal to the window size results in \SI{0}{\percent} overlap, while a hop size equal to half the window size results in \SI{50}{\percent} overlap.
\par
Each segment is multiplied by a window function that tapers the signal to zero at the boundaries, thus satisfying the periodicity requirement of the FFT and emphasizing the center of the segment. Zero-padding is then applied by appending additional zeros to each segment to increase the frequency resolution; the amount of padding is controlled by the zero-padding factor. The FFT is then applied to each padded segment, producing a complex-valued frequency vector. These vectors are stored as columns in a matrix, resulting in a matrix with as many columns as there are segments.
\par
The STFT output is generally not trivially invertible and depends significantly on the chosen parameters, reflecting a trade-off between time and frequency resolution. Therefore, the output of each experiment run is this matrix of complex vectors~\cite[src/Experiments.jl, run\_stft()]{ngafftb}, and the corresponding reference solution is computed using \texttt{float128} arithmetic. For this study, we adopt standard parameters tailored to a \SI{44.1}{\kilo\hertz} sampling rate: a window size of 2048 samples, a hop size of 1024 samples (corresponding to \SI{50}{\percent} overlap), the \textsc{Hann} window function and a zero-padding factor of two (yielding segment lengths twice the window size).
\section{Results}\label{sec:results}
In this section, we present the results of the experiments described above. For each experiment, the computed relative errors are sorted and displayed as cumulative error distributions. This visualization enables a straightforward assessment of the proportion of runs that remained below a given relative error threshold, while also supporting meaningful comparisons across the diverse datasets.
\subsection{Heat Equation Experiments}
\begin{figure*}[tbp]
	\begin{center}
        \subfloat[$N_x = 10^2$]{
			\includegraphics{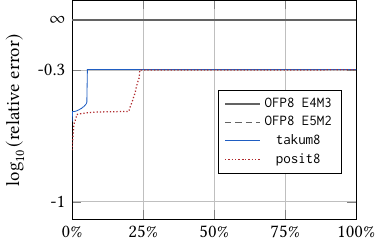}
			\includegraphics{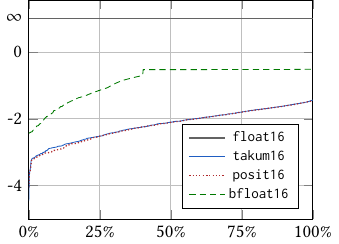}
			\includegraphics{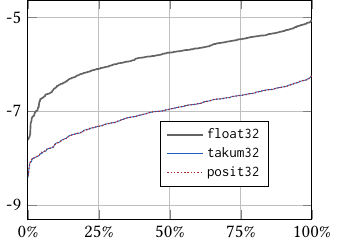}
        }\\
        \subfloat[$N_x = 10^3$]{
			\includegraphics{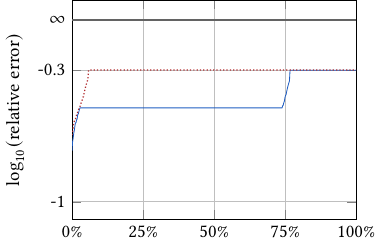}
			\includegraphics{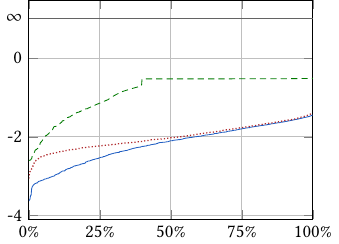}
			\includegraphics{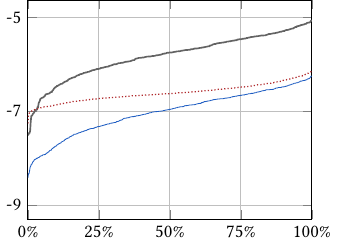}
        }\\
        \subfloat[$N_x = 5 \cdot 10^3$]{
			\includegraphics{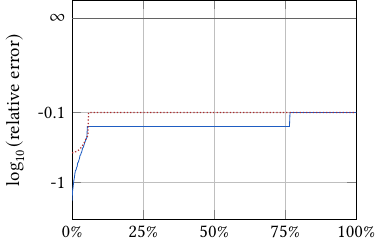}
			\includegraphics{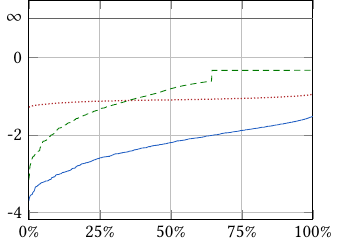}
			\includegraphics{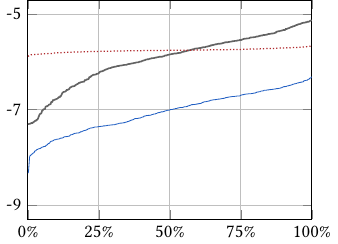}
        }\\
        \subfloat[$N_x = 10^4$]{
			\includegraphics{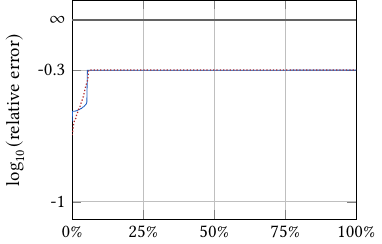}
			\includegraphics{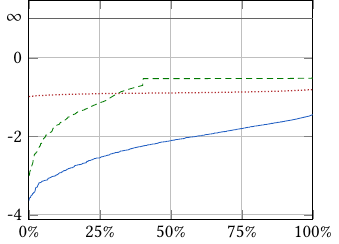}
			\includegraphics{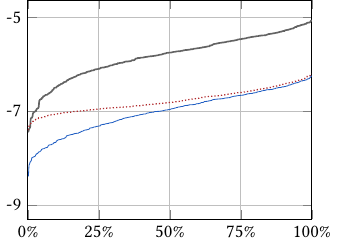}
        }
	\end{center}
	\caption{
		Cumulative error distribution of the relative errors of
		the heat equation PDE test problem solutions computed
		using a range of machine number types, with the specified
		space discretization parameter $N_x$ and $N_t \in \{ 1,\dots,500 \}$.
		The symbol $\infty$ denotes an overflow in the underlying number type.
	}
	\label{fig:pde-heat}
\end{figure*}
The results of the heat equation PDE test problem are presented in Figure~\ref{fig:pde-heat}, using four different values for the spatial discretization parameter $N_x$. For each value of $N_x$, the temporal discretization parameter $N_t$---which primarily influences numerical stability---was varied between $1$ and $500$. Each column in the figure corresponds to a numerical type with 8-, 16-, or 32-bit precision. The 64-bit types are omitted, as their numerical error consistently fell below the discretization error, resulting in identical outputs across all test cases.
\par
For the 8-bit formats, it is evident that OFP8 consistently overflows and fails to produce any meaningful result. In contrast, both posit and takum formats yield finite outputs across all configurations. At $N_x = 10^2$, posits outperform takums; however, for $N_x \in \{10^3, 5 \cdot 10^3\}$, takums surpass posits, while for $N_x = 10^4$, the two formats perform comparably.
\par
At 16-bit precision, \texttt{float16} also suffers from consistent overflow, yielding no useful results. In contrast, \texttt{bfloat16} produces finite outputs. For $N_x = 10^2$, posits and takums perform similarly and achieve at least an order of magnitude better accuracy than \texttt{bfloat16}. In all other cases, takums outperform both posits and \texttt{bfloat16}, while posits underperform relative to \texttt{bfloat16} for $N_x \in \{5 \cdot 10^3, 10^4\}$.
\par
At 32-bit precision, both \texttt{float32} and \texttt{takum32} deliver consistent and stable results across all scenarios, as expected due to the dominance of discretization error at this resolution. However, \texttt{posit32} exhibits erratic behavior depending on the value of $N_x$, ranging from matching \texttt{takum32} in accuracy to exceeding both \texttt{float32} and \texttt{takum32}. This variability may be attributed to the posit format's reduced precision for values with very large or very small magnitudes.
\subsection{Poisson Equation Experiments}
\begin{figure*}[tbp]
	\begin{center}
        \subfloat[$\sigma = 0.1$]{
			\includegraphics{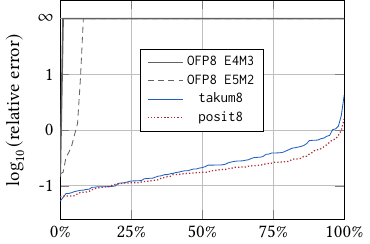}
			\includegraphics{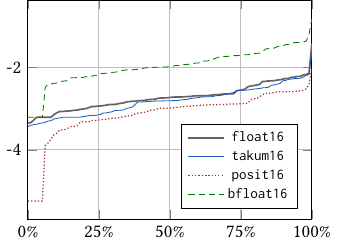}
			\includegraphics{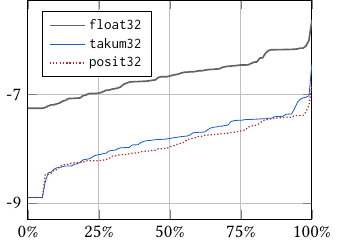}
        }\\
        \subfloat[$\sigma = 0.2$]{
			\includegraphics{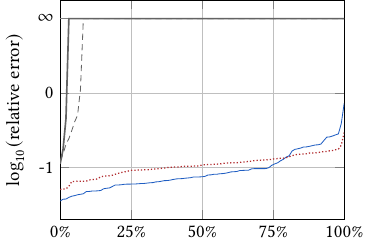}
			\includegraphics{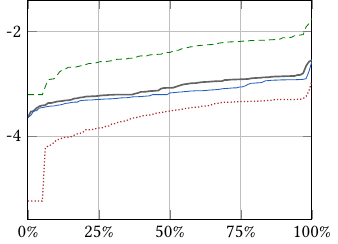}
			\includegraphics{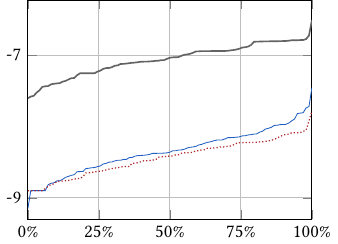}
        }\\
        \subfloat[$\sigma = 0.3$]{
			\includegraphics{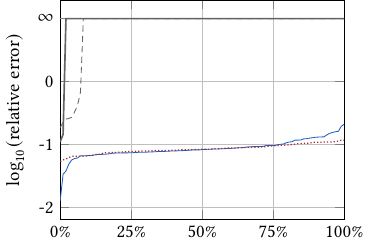}
			\includegraphics{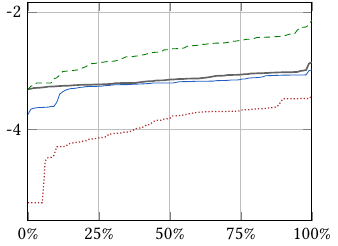}
			\includegraphics{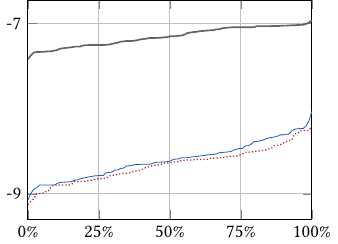}
        }\\
        \subfloat[$\sigma = 0.4$]{
			\includegraphics{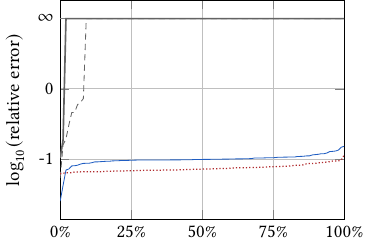}
			\includegraphics{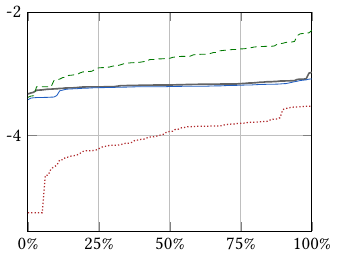}
			\includegraphics{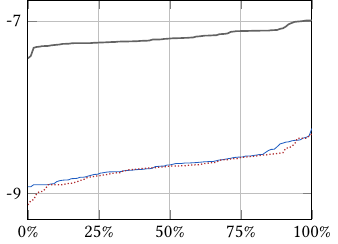}
        }
	\end{center}
	\caption{
		Cumulative error distribution of the relative errors of
		the \textsc{Poisson} equation PDE test problem solutions computed
		using a range of machine number types, with the specified
		right-hand side parameter $\sigma$ and $N_x \in \{ 2,\dots,100 \}$.
		The symbol $\infty$ denotes an overflow in the underlying number type.
	}
	\label{fig:pde-poisson}
\end{figure*}
Figure~\ref{fig:pde-poisson} presents the results of the \textsc{Poisson} equation PDE test problem. In this case, the parameter $\sigma$ was varied across four values, and for each, the space discretization parameter was swept from $2$ to $100$. Each column corresponds to a number type with 8-, 16-, or 32-bit precision. As with the heat equation results, 64-bit types are omitted since their numerical errors fell below the discretization error threshold in all cases.
\par
At 8-bit precision, the OFP8 types consistently failed to produce finite results, with the exception of a few cases for E5M2. In contrast, both posits and takums yielded finite outputs with reasonable relative errors. Their overall performance is similar, with takums slightly outperforming posits at $\sigma = 0.2$, and posits slightly ahead at $\sigma = 0.4$.
\par
At 16-bit precision, \texttt{bfloat16} underperforms relative to \texttt{float16}, whereas \texttt{takum16} slightly outperforms it. Notably, \texttt{posit16} achieves significantly better accuracy than all other types at this precision level.
\par
At 32-bit precision, posits and takums again perform comparably, both achieving one to two orders of magnitude lower relative error than \texttt{float32}.
\subsection{Image FFT Experiments}
\begin{figure*}[tbp]
	\begin{center}
        \subfloat[8 bits]{
			\includegraphics{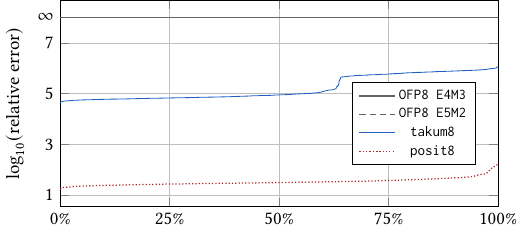}
    	}
    	\subfloat[16 bits]{
			\includegraphics{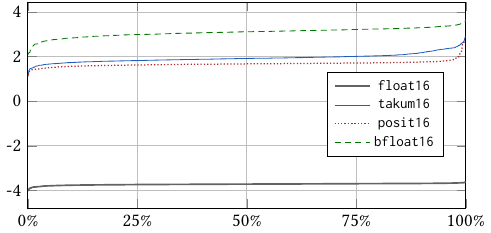}
        }\\
        \subfloat[32 bits]{
        	\includegraphics{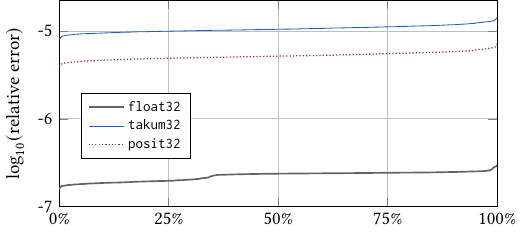}
        }
        \subfloat[64 bits]{
			\includegraphics{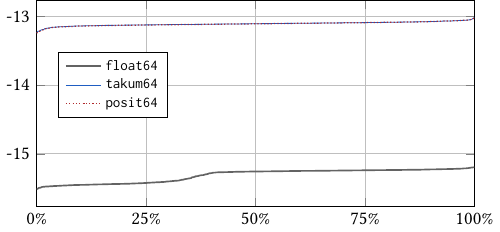}
        }
	\end{center}
	\caption{
		Cumulative error distribution of the relative errors of
		the image FFT solutions computed
		using a range of machine number types.
		The symbol $\infty$ denotes an overflow in the underlying number type.
	}
	\label{fig:images}
\end{figure*}
Figure~\ref{fig:images} shows the results for the image dataset after the round-trip FFT conversions, ranging from 8 to 64 bits. These results confirm that the OFP8 types are unsuitable for this task. At 16 bits, \texttt{float16} exhibits clear superiority, while posits and takums perform similarly and outperform \texttt{bfloat16}. At 32 bits, \texttt{float32} leads in performance, with posits and takums lagging by approximately an order of magnitude. This trend continues at 64 bits, where posits and takums perform similarly, but both fall approximately two orders of magnitude behind \texttt{float64}.
\subsection{Audio STFT Experiments}
\begin{figure*}[tbp]
	\begin{center}
        \subfloat[8 bits]{
			\includegraphics{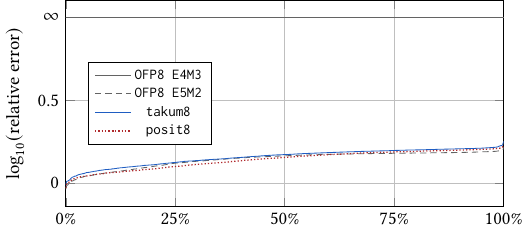}
    	}
    	\subfloat[16 bits]{
			\includegraphics{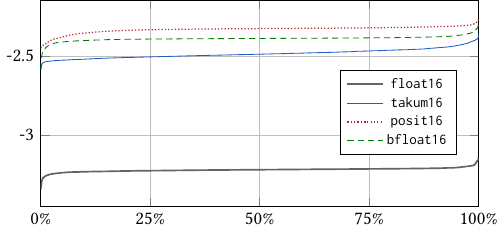}
        }\\
        \subfloat[32 bits]{
			\includegraphics{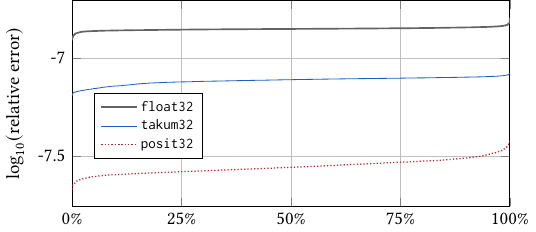}
        }
        \subfloat[64 bits]{
			\includegraphics{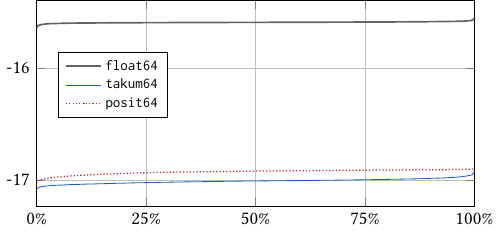}
        }
	\end{center}
	\caption{
		Cumulative error distribution of the relative errors of
		the audio STFT solutions computed
		using a range of machine number types.
		The symbol $\infty$ denotes an overflow in the underlying number type.
	}
	\label{fig:audios}
\end{figure*}
The results for the audio dataset, shown in Figure~\ref{fig:audios}, again confirm that OFP8 is unsuitable for this type of application. In contrast, \texttt{posit8} and \texttt{takum8} exhibit identical performance. At 16 bits, \texttt{float16} is nearly an order of magnitude better than the other types, which are dominated by takums, followed by \texttt{bfloat16} and posits, which perform the worst. At 32 bits, all types perform very similarly, with posits leading, followed closely by takums, and \texttt{float32} being the least effective. At 64 bits, posits and takums perform almost identically, yielding results over an order of magnitude better than \texttt{float64}.
\section{Conclusion}\label{sec:conclusion}
In this paper, we evaluated IEEE 754 floating-point numbers, OFP8 E4M3 and E5M2, \texttt{bfloat16}, posits, and takums for use in spectral methods for solving partial differential equations (PDEs) and signal analysis. Our experiments demonstrate that the OFP8 types are unsuitable for this type of application. Additionally, \texttt{bfloat16} outperforms \texttt{float16} due to its ability to avoid overflow, despite some numerical errors. The tapered-precision formats, posits and takums, generally outperform OFP8 and \texttt{bfloat16}, with takums consistently performing better than posits, which occasionally fall behind \texttt{bfloat16}. This suggests that takums are likely the best choice at 16 bits, while at 8 bits, the results depend on the specific problem, with posits occasionally performing better. At 32 and 64 bits, the results are more mixed, likely due to the approaching intrinsic discretization error. However, overall, posits and takums---except in a few cases---outperform \texttt{float32} and \texttt{float64}, respectively.
\par
These findings are significant, as they suggest that \texttt{takum16} is a viable candidate for replacing \texttt{bfloat16} as the state-of-the-art in 16-bit arithmetic. They also indicate that the high dynamic range of takums does not significantly hinder their numerical performance; rather, it appears to enhance stability. Posits were also shown to provide overall superior performance with a few exceptions.
\par
Future research could extend the benchmark methods to include additional PDE types, diverse initial conditions and right-hand sides, varied datasets, and different STFT parameters.
\section*{Author Contributions}
Specific roles and contributions will be detailed following the review process.
\printbibliography
\end{document}